\documentclass{amsart}

\usepackage{amssymb,amsmath,amsthm,latexsym,stmaryrd}

\newtheorem{theorem}{Theorem}[section]
\newtheorem*{theorem A}{Theorem A}
\newtheorem*{theorem B}{N\"olker's Theorem}

\theoremstyle{remark}

\theoremstyle{remark}

\theoremstyle{definition}

\numberwithin{equation}{section}
\def\({\left ( }
    \def\){\right )}
\def\<{\left < }
\def\>{\right >}

\newcommand{\ie}{i.e. }
\newcommand{\f}{\varphi}

\newcommand{\g}{\mathfrak{g}}
\newcommand{\n}{\nabla}

\newcommand{\R}{\mathbb R}

\newcommand{\F}{\mathcal{F}}

\newcommand{\ta}{\theta}
\newcommand{\om}{\omega}
\newcommand{\lm}{\lambda}

\newcommand{\al}{\alpha}
\newcommand{\bt}{\beta}
\newcommand{\Dt}{\Delta}

\newcommand{\Ker}{\mathrm{ker}}

\newcommand{\Span}{\operatorname{span}}
\newcommand{\Id}{\operatorname{Id}}
\newcommand{\tr}{\operatorname{tr}}
\newcommand{\de}{\operatorname{det}}

\newcommand{\RNum}[1]{\uppercase\expandafter{\romannumeral #1\relax}}
\newcommand{\Bia}[1]{\mathrm{Bia}(\mathrm{\uppercase\expandafter{\romannumeral #1\relax}})}

\begin{document}

%

\vspace{2cm}

\title{Matrix Lie Groups as 3-Dimensional Almost Contact B-Metric Manifolds}

\author{Hristo Manev$^{1,2}$}
\address[1]{Medical University of Plovdiv, Faculty of Pharmacy,
Department of Pharmaceutical Sciences,   15-A Vasil Aprilov
Blvd.,   Plovdiv 4002,   Bulgaria;}  \address[2]{Paisii
Hilendarski University of Plovdiv,   Faculty of Mathematics and
Informatics,   Department of Algebra and Geometry,   236
Bulgaria Blvd.,   Plovdiv 4027,   Bulgaria}
\email{hmanev@uni-plovdiv.bg}



\dedicatory{{\rm Dedicated to my doctoral advisor Prof. Dimitar Mekerov on the occasion of his 70th birthday}}


\begin{abstract}
The object of investigation are Lie groups considered as almost contact B-metric manifolds of the lowest dimension three.
It is established a correspondence of all basic-class-manifolds of the Ganchev-Mihova-Gribachev classification of the studied manifolds and the explicit matrix representation of Lie groups. Some known Lie groups are equiped with almost contact B-metric structure of different types.
\end{abstract}
\maketitle
\section*{Introduction}
The differential geometry of almost contact manifolds
is well studied (e.g. \cite{Blair}). In \cite{GaMiGr}, the almost contact manifolds with B-metric are introduced and classified. These manifolds can be studied as the odd-dimensional
counterpart of the almost complex manifolds with Norden metric \cite{GaBo, GrMeDj}.

As it is known from \cite{Gil}, every representation of a Lie algebra corresponds uniquely to a representation of the simply connected Lie group. This relation is one-to-one. Then, knowing the representation of a certain Lie algebra settles the question of the representation of its Lie group.

In \cite{IcoMek}, an arbitrary Lie group considered as a manifold is equiped with an almost contact B-metric structure. For the lowest dimension three, it is deduced an equivalence between the classification of considered manifolds given in \cite{GaMiGr} and the corresponding Lie algebra determined by commutators.

In the present work, it is found the explicit correspondence between all basic-class-manifolds of the classification in \cite{GaMiGr} and Lie groups given by their explicit matrix representation.

The paper is organized as follows. In
Sect.~\ref{sect-prel} we recall some preliminary facts about the almost contact B-metric manifolds  and  Lie algebras related to them. In Sect.~\ref{sect-lie} we give the correspondence of all basic-class-manifolds and matrix Lie groups.
Sect.~\ref{sect-ex} is devoted to some examples connected to the previous investigations.

\section{Preliminaries}\label{sect-prel}

\subsection{Almost contact B-metric manifolds}

Let us denote an almost contact B-met\-ric manifold by $(M,\f,\xi,\eta,g)$, \ie $M$
is an $(2n+1)$-dimensional differentiable manifold, $\f$ is an endomorphism
 of the tangent bundle, $\xi$ is a Reeb vector field, $\eta$ is its dual contact 1-form and  $g$ is a pseu\-do-Rie\-mannian
metric (called a \emph{B-metric}) of signature $(n+1,n)$ such that \cite{GaMiGr}
\[
\begin{array}{c}
\f\xi = 0,\qquad \f^2 = -\Id + \eta \otimes \xi,\qquad
\eta\circ\f=0,\qquad \eta(\xi)=1,\\[0pt]
g(\f x, \f y) = - g(x,y) + \eta(x)\eta(y),
\end{array}
\]
where $\Id$ is the identity. In the latter equalities and further, $x$, $y$, $z$, $w$ will stand for arbitrary elements of the
algebra of the smooth vector fields on $M$ or vectors in the tangent space $T_pM$ of $M$ at an arbitrary
point $p$ in $M$.

In \cite{GaMiGr}, it is given a classification of almost contact B-metric manifolds, consisting of eleven basic
classes $\F_1$, $\F_2$, $\dots$, $\F_{11}$. It is made with respect
to the tensor $F$ of type (0,3) defined by
\begin{equation*}\label{F=nfi}
F(x,y,z)=g\bigl( \left( \nabla_x \f \right)y,z\bigr),
\end{equation*}
where $\n$ stands for the Levi-Civita connection
of $g$.

The special class $\F_0$, determined by the condition $F(x,y,z)=0$, is the intersection of all basic classes and it is known as the class of the \emph{cosymplectic B-metric manifolds}.

Let $\left\{\xi;e_i\right\}$ $(i=1,2,\dots,2n)$ be a basis of
$T_pM$ and let $\left(g^{ij}\right)$ be the inverse matrix of
$\left(g_{ij}\right)$. Then the Lee 1-forms
$\theta$, $\theta^*$, $\omega$ associated with $F$ are  defined
by
\[
\theta(z)=g^{ij}F(e_i,e_j,z),\quad \theta^*(z)=g^{ij}F(e_i,\f
e_j,z), \quad \omega(z)=F(\xi,\xi,z).
\]

In this work, we consider the case of the lowest dimension of the considered manifolds, \ie $\dim{M}=3$.

We introduce an almost contact structure $(\f,\xi,\eta)$ on $M$ defined by
\begin{equation*}\label{strL}
\begin{array}{c}
\f e_1=e_{2},\quad \f e_{2}=- e_1,\quad \f e_0=0,\quad \xi=
e_0,\quad \\[0pt]
\eta(e_1)=\eta(e_{2})=0,\quad \eta(e_0)=1
\end{array}
\end{equation*}
and a B-metric $g$ such that
\begin{equation*}\label{gij}
g(e_0,e_0)=g(e_1,e_1)=-g(e_2,e_2)=1,\quad g(e_i,e_j)=0,\;\; i\neq j
\in \{0,1,2\}.
\end{equation*}

The components of $F$, $\ta$, $\ta^*$, $\om$ with respect to the \emph{$\f$-basis}
$\left\{e_0,e_1,e_2\right\}$ are denoted by $F_{ijk}=F(e_i,e_j,e_k)$, $\ta_k=\ta(e_k)$, $\ta^*_k=\ta^*(e_k)$, $\om_k=\om(e_k)$.
According to \cite{HM}, we have:
\[
\begin{array}{lll}
\ta_0=F_{110}-F_{220},\qquad &\ta_1=F_{111}-F_{221},\qquad &\ta_2=F_{112}-F_{211},\\[0pt]
\ta^*_0=F_{120}+F_{210},\qquad &\ta^*_1=F_{112}+F_{211},\qquad &\ta^*_2=F_{111}+F_{221},\\[0pt]
\om_0=0,\qquad &\om_1=F_{001},\qquad &\om_2=F_{002}.
\end{array}
\]

If $F_s$ $(s=1,2,\dots,11)$ are the components of $F$ in the
corresponding basic classes $\F_s$ then: \cite{HM}
\begin{equation}\label{Fi3}
\begin{array}{l}
F_{1}(x,y,z)=\left(x^1\ta_1-x^2\ta_2\right)\left(y^1z^1+y^2z^2\right),\\[0pt]
\qquad\ta_1=F_{111}=F_{122},\qquad \ta_2=-F_{211}=-F_{222}; \\[0pt]
F_{2}(x,y,z)=F_{3}(x,y,z)=0;
\\[0pt]
F_{4}(x,y,z)=\frac{1}{2}\ta_0\Bigl\{x^1\left(y^0z^1+y^1z^0\right)
-x^2\left(y^0z^2+y^2z^0\right)\bigr\},\\[0pt]
\qquad \frac{1}{2}\ta_0=F_{101}=F_{110}=-F_{202}=-F_{220};\\[0pt]
F_{5}(x,y,z)=\frac{1}{2}\ta^*_0\bigl\{x^1\left(y^0z^2+y^2z^0\right)
+x^2\left(y^0z^1+y^1z^0\right)\bigr\},\\[0pt]
\qquad \frac{1}{2}\ta^*_0=F_{102}=F_{120}=F_{201}=F_{210};\\[0pt]
F_{6}(x,y,z)=F_{7}(x,y,z)=0;\\[0pt]
F_{8}(x,y,z)=\lm\bigl\{x^1\left(y^0z^1+y^1z^0\right)
+x^2\left(y^0z^2+y^2z^0\right)\bigr\},\\[0pt]
\qquad \lm=F_{101}=F_{110}=F_{202}=F_{220};\\[0pt]
F_{9}(x,y,z)=\mu\bigl\{x^1\left(y^0z^2+y^2z^0\right)
-x^2\left(y^0z^1+y^1z^0\right)\bigr\},\\[0pt]
\qquad \mu=F_{102}=F_{120}=-F_{201}=-F_{210};\\[0pt]
F_{10}(x,y,z)=\nu x^0\left(y^1z^1+y^2z^2\right),\qquad
\nu=F_{011}=F_{022};\\[0pt]
F_{11}(x,y,z)=x^0\bigl\{\left(y^1z^0+y^0z^1\right)\om_{1}
+\left(y^2z^0+y^0z^2\right)\om_{2}\bigr\},\\[0pt]
\qquad \om_1=F_{010}=F_{001},\qquad \om_2=F_{020}=F_{002},
\end{array}
\end{equation}
where $x=x^ie_i$, $y=y^je_j$, $z=z^ke_k$.
Obviously, the class of three-dimensional almost contact B-metric
manifolds is
\[
\F_1 \oplus \F_4 \oplus \F_5 \oplus \F_8 \oplus \F_9 \oplus
\F_{10} \oplus \F_{11}.
\]

\subsection{The Lie algebras corresponding to Lie groups with  almost contact B-metric structure }

Let $L$ and $\mathfrak{l}$ be a 3-dimensional real connected Lie group and its corresponding Lie algebra.
If $\{E_{0},E_{1},E_{2}\}$ is a basis of left invariant vector fields then an almost contact structure $(\f,\xi,\eta)$ and a B-metric $g$ are defined as follows:
\begin{equation*}\label{strLgL}
\begin{array}{l}
\f E_0=0,\quad \f E_1=E_{2},\quad \f E_{2}=- E_1,\quad \xi=
E_0,\quad \\[4pt]
\eta(E_0)=1,\quad \eta(E_1)=\eta(E_{2})=0,\quad \\[4pt]
g(E_0,E_0)=g(E_1,E_1)=-g(E_{2},E_{2})=1, \\[4pt]
g(E_0,E_1)=g(E_0,E_2)=g(E_1,E_2)=0.
\end{array}
\end{equation*}
Then $(L,\f,\xi,\eta,g)$ is a 3-dimensional almost contact B-metric manifold.
If its Lie algebra $\mathfrak{l}$ is determined by
\begin{equation*}\label{lie}
\left[E_{i},E_{j}\right]=C_{ij}^k E_{k}, \quad i, j, k \in \{0,1,2\},
\end{equation*}
then the components of $F$, $\ta$, $\ta^*$, $\om$ are the following: \cite{IcoMek}
\begin{equation*}\label{FijkC}
\begin{array}{l}
\begin{array}{ll}
F_{111}=F_{122}=2C_{12}^1,\quad &
F_{211}=F_{222}=2C_{12}^2,\\[4pt]
F_{120}=F_{102}=-C_{01}^1,\quad &
F_{020}=F_{002}=-C_{01}^0,\\[4pt]
F_{210}=F_{201}=-C_{02}^2,\quad &
F_{010}=F_{001}=C_{02}^0,\\[4pt]
\end{array}\ \\[4pt]
\begin{array}{l}
F_{110}=F_{101}=\frac12 \left(C_{12}^0-C_{01}^2+C_{02}^1\right),\\[4pt]
F_{220}=F_{202}=\frac12 \left(C_{12}^0+C_{01}^2-C_{02}^1\right),\\[4pt]
F_{011}=F_{022}=C_{12}^0+C_{01}^2+C_{02}^1,
\end{array}
\end{array}
\end{equation*}
\begin{equation*}\label{titiC}
\begin{array}{l}
\begin{array}{l}
\ta_{0}=-C_{01}^2+C_{02}^1,\\[4pt]
\ta_{1}=2C_{12}^1,\\[4pt]
\ta_{2}=-2C_{12}^2,\\[4pt]
\end{array}\quad
\begin{array}{l}
\ta^*_{0}=-C_{01}^1-C_{02}^2,\\[4pt]
\ta^*_{1}=2C_{12}^2,\\[4pt]
\ta^*_{2}=2C_{12}^1,\\[4pt]
\end{array}\quad
\begin{array}{l}
\om_{0}=0,\\[4pt]
\om_{1}=C_{02}^0,\\[4pt]
\om_{2}=-C_{01}^0.
\end{array}
\end{array}
\end{equation*}

\begin{theorem A}[\cite{IcoMek}]\label{thm-Fi-L}
The manifold $(L,\f,\xi,\eta,g)$ belongs to the basic class $\F_s$
($s \in \{1,4,5,8,9,10,11\}$) if and only
if the corresponding Lie algebra $\mathfrak{l}$ is determined by
the following commutators:
\begin{equation*}\label{Fi-L}
\begin{array}{ll}
\F_1:\; &[E_0,E_1]=[E_0,E_2]=0, \quad [E_1,E_2]=\al E_1+\bt E_2;
\\[4pt]
\F_4:\; &[E_0,E_1]=\al E_2, \quad [E_0,E_2]=-\al E_1, \quad
[E_1,E_2]=0;
\\[4pt]
\F_5:\; &[E_0,E_1]=\al E_1, \quad [E_0,E_2]=\al E_2, \quad
[E_1,E_2]=0;
\\[4pt]
\F_8:\; &[E_0,E_1]=\al E_2, \quad [E_0,E_2]=\al E_1, \quad
[E_1,E_2]=-2\al E_0;
\\[4pt]
\F_9:\; &[E_0,E_1]=\al E_1, \quad [E_0,E_2]=-\al E_2, \quad
[E_1,E_2]=0;
\\[4pt]
\F_{10}:\; &[E_0,E_1]=\al E_2, \quad [E_0,E_2]=\al E_1, \quad
[E_1,E_2]=0;
\\[4pt]
\F_{11}:\; &[E_0,E_1]=\al E_0, \quad [E_0,E_2]=\bt E_0, \quad
[E_1,E_2]=0,
\end{array}
\end{equation*}
where $\al$, $\bt$ are arbitrary real parameters. Moreover, the
relations of $\al$ and $\bt$ with the non-zero components
$F_{ijk}$ in the different basic classes $\F_s$ from \eqref{Fi3} are the
following:
\begin{equation*}\label{Fi-L-alpha}
\begin{array}{ll}
\F_1:\quad \al=\frac12 \ta_1,\quad \bt=\frac12 \ta_2; \qquad
&\F_4:\quad \al=\frac12 \ta_0;\\[4pt]
\F_5:\quad \al=-\frac12 \ta^*_0; \qquad
&\F_8:\quad \al=-\lm; \\[4pt]
\F_9:\quad \al=-\mu; \qquad
&\F_{10}:\quad \al=\frac12 \nu;   \\[4pt]
\F_{11}:\quad \al=-\om_2, \quad \bt=\om_1.\qquad &
\end{array}
\end{equation*}
\end{theorem A}

\section{The main result}\label{sect-lie}

Theorem A yields an equivalence between the  manifolds of the classification in \cite{GaMiGr} and the corresponding Lie algebra.

It is known (e.g. \cite{Gil}) that for a real Lie algebra of finite dimension
there is a corresponding connected simply connected Lie group, which is determined uniquely up to isomoprhism.

It is arised the problem for determination of the Lie group which is isomorphic to the given Lie group $L$ equiped with a structure $(\f,\xi,\eta,g)$ in the class $\F_s$.

\begin{theorem}\label{thm:main}
Let $(L,\f,\xi,\eta,g)$ be an almost contact B-metric manifold belonging to the class $\F_s$
($s \in \{1,4,5,8,9,10,11\}$). Then the compact simply connected Lie group $G$ isomorphic to $L$, both with one and the same Lie algebra, has the form
\begin{equation*}\label{eA}
  e^A=E+tA+uA^2,
\end{equation*}
where $E$ is the identity matrix and $A$ is the matrix representation of the corresponding Lie algebra. The matrix form of $A$ as well as the real parameters $t$ and $u$  for the different classes $\F_s$ are given in Table~1, where $a,b,c\in\R$ and $\al$, $\bt$ are introduced in Theorem~A.
\end{theorem}

\begin{center}
\begin{table}
  \caption{The matrix form of $A$ and the values of $t$ and $u$ for $\F_s$}\label{T1}
\begin{tabular}{|l|l|l|}
\hline
$
\F_1:$ & $
        A=\left(
      \begin{array}{ccc}
        0 & 0 & 0 \\
        0 & \al b & \bt b \\
        0 & -\al a & -\bt a
      \end{array}
        \right)\quad $ &
$        t=
        \left\{
        \begin{array}{ll}
            \frac{e^{\tr{A}}-1}{\tr{A}}, & \tr{A}\neq 0 \\
                1, & \tr{A}= 0
        \end{array}
        \right.$\\
& $
        \tr{A}=\al b-\bt a $
        &
        $u=0
$
\\
\hline
$
\F_4:$&$
A=\left(
      \begin{array}{ccc}
        0 & -\al b & \al a \\
        0 & 0 & -\al c \\
        0 & \al c & 0 \\
      \end{array}
    \right)\quad $
    &
    $        t=
\left\{
  \begin{array}{ll}
    \frac{\sin{\sqrt{-\frac12 \tr{A^2}}}}{\sqrt{-\frac12 \tr{A^2}}}, & \tr{A^2}\neq 0 \\
    1, & \tr{A}= 0
  \end{array}
\right.$\\
& $
\tr{A^2}=-2\al^2 c^2 $ &
$    u=
\left\{
  \begin{array}{ll}
    \frac{1-\cos{\sqrt{-\frac12 \tr{A^2}}}}{-\frac12 \tr{A^2}}, & \tr{A^2}\neq 0 \\
    0, & \tr{A}= 0
  \end{array}
\right.
$
\\
\hline
$
\F_5:$&$
A=\left(
      \begin{array}{ccc}
        0 & \al a & \al b \\
        0 & -\al c & 0 \\
        0 & 0 & -\al c \\
      \end{array}
    \right)\quad $ & $
    t=
\left\{
  \begin{array}{ll}
    \frac{e^{\frac12 \tr{A}}-1}{\frac12 \tr{A}}, & \tr{A}\neq 0 \\
    1, & \tr{A}= 0
  \end{array}
\right.$\\
& $
\tr{A}=-2\al c $ &
$        u=0
$
\\
\hline
$
\F_8:$ & $
A=\left(
      \begin{array}{ccc}
        0 & \al b & \al a \\
        -2\al b & 0 & -\al c \\
        2\al a & -\al c & 0 \\
      \end{array}
    \right)\;
$ &
$    t=
\left\{
  \begin{array}{ll}
    -\frac{\sin{\al\sqrt{| \Dt|}}}{\al\sqrt{| \Dt|}}, & \tr{A^2}< 0 \\
    1,                                                 & \tr{A^2}= 0 \\
    \frac{\sinh{\al\sqrt{\Dt}}}{\al\sqrt{\Dt}}, & \tr{A^2}> 0 \\
  \end{array}
\right.$\\
& $
\begin{array}{l}
\tr{A^2}=2\al^2\Dt\\[4pt]
\Dt=2a^2-2b^2+c^2
\end{array}
$ & $
    u=
\left\{
  \begin{array}{ll}
    \frac{\cos{\al\sqrt{| \Dt|}}-1}{\al^2\Dt}, & \tr{A^2}< 0 \\
    \frac12,                                                 & \tr{A^2}= 0 \\
    \frac{\cosh{\al\sqrt{\Dt}}-1}{\al^2\Dt}, & \tr{A^2}> 0 \\
  \end{array}
\right.
$
\\
\hline
$
\F_9:$ & $
A=\left(
      \begin{array}{ccc}
        0 & \al a & -\al b \\
        0 & -\al c & 0 \\
        0 & 0 & \al c \\
      \end{array}
    \right)\quad $ & $
    t=
\left\{
  \begin{array}{ll}
    \frac{\sinh{\sqrt{\frac12 \tr{A^2}}}}{\sqrt{\frac12 \tr{A^2}}}, & \tr{A^2}\neq 0 \\
    1, & \tr{A^2}= 0
  \end{array}
\right.$\\
& $
\tr{A^2}=2\al^2 c^2 $ & $
    u=
\left\{
  \begin{array}{ll}
    \frac{\cosh{\frac12 \tr{A^2}}-1}{\frac12 \tr{A^2}}, & \tr{A^2}\neq 0 \\
    0,                                                  & \tr{A^2}= 0 \\
  \end{array}
\right.
$
\\
\hline
$
\F_{10}:$ & $
A=\left(
      \begin{array}{ccc}
        0 & \al b & \al a \\
        0 & 0 & -\al c \\
        0 & -\al c & 0 \\
      \end{array}
    \right)\quad $ & $
    t=
\left\{
  \begin{array}{ll}
    \frac{\sinh{\al c}}{\al c}, & \tr{A^2}\neq 0 \\
    1,                        & \tr{A^2}= 0 \\
  \end{array}
\right.$\\
& $
\tr{A^2}=2\al^2 c^2 $ & $
    u=
\left\{
  \begin{array}{ll}
    \frac{\cosh{\al c}}{\al^2c^2}, & \tr{A^2}\neq 0 \\
    0,                         & \tr{A^2}= 0 \\
  \end{array}
\right.
$
\\
\hline
$
\F_{11}:$ & $
A=\left(
      \begin{array}{ccc}
        \al a + \bt b & 0 & 0 \\
        -\al c & 0 & 0 \\
        -\bt c & 0 & 0 \\
      \end{array}
    \right)\quad $ &
$    t=
\left\{
  \begin{array}{ll}
    \frac{e^{-\tr{A}}-1}{\tr{A}}, & \tr{A}\neq 0 \\
    1,                        & \tr{A}= 0 \\
  \end{array}
\right.$\\
& $
\tr{A}=\al a+\bt b $ & $
    u=0
$\\
\hline
\end{tabular}
\end{table}
\end{center}

\subsection{Proof of the theorem}

\subsubsection{Lie groups as manifolds from the classes $\F_1$, $\F_5$, $\F_{11}$}\label{2.1}

Firstly, let us consider the case when $(L,\f,\xi,\eta,g)$ is a $\F_1$-manifold.
Then, according to Theorem A, the corresponding Lie algebra is determined by

\begin{equation}\label{com1}
[E_0,E_1]=[E_0,E_2]=0,\quad [E_1,E_2]=\al E_1+\bt E_2,
\end{equation}
where $\al=\frac12 \ta_1$, $\bt=\frac12 \ta_2$.

From \eqref{com1} we have the nonzero values of the commutation coefficients:
\begin{equation}\label{Cij}
C_{12}^1=-C_{21}^1=\al,\quad C_{12}^2=-C_{21}^2=-\bt.
\end{equation}

According to \cite{Gil}, the commutation coefficients provide a matrix representation of the Lie algebra.
This representation is obtained by the basic matrices $M_i$, which entries are determined by
\begin{equation}\label{Mij}
(M_i)_j^k=-C_{ij}^k.
\end{equation}
Using \eqref{Cij} and \eqref{Mij}, we obtain

\begin{equation*}\label{}
M_0=\left(
      \begin{array}{ccc}
        0 & 0 & 0 \\
        0 & 0 & 0 \\
        0 & 0 & 0 \\
      \end{array}
    \right),\quad
M_1=\left(
      \begin{array}{ccc}
        0 & 0 & 0 \\
        0 & 0 & 0 \\
        0 & -\al & -\bt\\
      \end{array}
    \right),\quad
M_2=\left(
      \begin{array}{ccc}
        0 & 0 & 0 \\
        0 & \al & \bt \\
        0 & 0 & 0 \\
      \end{array}
    \right).
\end{equation*}

Let us suppose that $(a,b)\neq (0,0)$.
The matrix representation of the considered Lie algebra $\mathfrak{g}_1$ is the matrix $A$, given in Table~1.
Then, the characteristic polynomial of $A$ is
\begin{equation*}\label{}
P_A(\lm)= \lm^2(\lm - \al b + \bt a).
\end{equation*}
Therefore, the eigenvalues $\lm_i$ ($i={1,2,3}$) are
\begin{equation*}\label{}
\lm_1=\lm_2=0, \qquad \lm_3=\al b - \bt a.
\end{equation*}
Hence, the corresponding linearly independent eigenvectors $p_i$ ($i={1,2,3}$) are
\begin{equation*}\label{}
p_1(1,0,0)^T, \qquad p_2(0,\bt,-\al)^T, \qquad p_3(0,-b,a)^T
\end{equation*}
and their matrix $P$ has the following form
\begin{equation*}\label{}
P=\left(
      \begin{array}{ccc}
        1 & 0 & 0 \\
        0 & \bt & -b \\
        0 & -\al & a
      \end{array}
        \right).
\end{equation*}
Let us denote  $\Dt=\bt a-\al b $. Then, we have $\Dt=\de P=-\tr A$.

Let us consider the first case $\Dt \neq 0$.
Then we have
\begin{equation*}\label{}
P^{-1}=\left(
      \begin{array}{ccc}
        1 & 0 & 0 \\
        0 & \frac{a}{\Dt} & \frac{b}{\Dt} \\
        0 & \frac{\al}{\Dt} & \frac{\bt}{\Dt}
      \end{array}
        \right).
\end{equation*}
The Jordan matrix $J$ is the diagonal matrix with elements $J_{ii}=\lm_i$.
It is known that
$e^A=Pe^JP^{-1}$.
Then we obtain the matrix Lie group representation $G_1$ of the considered Lie algebra $\mathfrak{g}_1$ in this case as follows
\begin{equation*}\label{}
G_1=\left\{
e^A=\left(\left.
      \begin{array}{ccc}
        1 & 0 & 0 \\
        0 & 1+\frac{\al b}{\Dt}\left(1-e^{-\Dt}\right)
          & \frac{\bt b}{\Dt}\left(1-e^{-\Dt}\right) \\
        0 & -\frac{\al a}{\Dt}\left(1-e^{-\Dt}\right)
          & 1-\frac{\bt a}{\Dt}\left(1-e^{-\Dt}\right) \\
      \end{array}
    \right)\;\right| \;
      \begin{array}{l}
                \Dt\neq 0
      \end{array}
    \right\},
\end{equation*}
which can be written as
\begin{equation*}\label{}
G_1=\left\{\left.
e^A=E+\frac{1-e^{-\Dt}}{\Dt}A\;\right| \;
     \Dt\neq 0
    \right\}.
\end{equation*}

In the second case $\Dt  = 0$, the matrix $P$ is non-invertible.
Then $A$ is nilpotent and $e^A$ can be computed directly from
\begin{equation*}\label{}
e^A=E+A+\frac{A^2}{2!}+\frac{A^3}{3!}+\dots++\frac{A^{q-1}}{(q-1)!},
\end{equation*}
where $q$ stands for the degree of $A$.
Then, from the form of $A$ in Table~1 we obtain that $q=2$.
Therefore, the matrix representation $G_1$ of the Lie group for the Lie algebra $\mathfrak{g}_1$ in this case is
\begin{equation*}\label{}
G_1=\left\{
e^A=E+A, \;
    \Dt= 0
    \right\},
\end{equation*}
\ie
\begin{equation*}\label{}
G_1=\left\{
e^A=\left(
      \begin{array}{ccc}
        1 & 0 & 0 \\
        0 & 1+\al b
          & \bt b \\
        0 & -\al a
          & 1-\bt a \\
      \end{array}
    \right), \;
    \Dt=0
    \right\}.
\end{equation*}

Both the cases for $\F_1$ can be generalized as it is shown in Table~1.

By similiar considerations we obtain the results in Table~1 for the classes $\F_5$ and $\F_{11}$.

\subsubsection{Lie groups as manifolds from the classes $\F_4$, $\F_9$, $\F_{10}$}Firstly, let us consider the case when $(L,\f,\xi,\eta,g)$ is a $\F_4$-manifold.
Then, according to Theorem A, the corresponding Lie algebra is determined by
\begin{equation}\label{C4}
[E_0,E_1]=\al E_2,\quad [E_0,E_2]=-\al E_1,\quad [E_1,E_2]=0,
\end{equation}
where $\al=\frac12 \ta_0$.

By similar way as in \S2.1.1, we obtain the form of $A$ for the corresponding Lie algebra $\g_4$, given in Table~1.
In this class  $\tr{A}=0$ and we have two cases according to $\tr{A^2}=-2\al^2 c^2$.
Then, the matrix representation is respectively
\begin{equation*}\label{}
G_4=\left\{
e^A=\left(
      \begin{array}{ccc}
        1 & \frac{a}{c}(1-\cos{\al c})-\frac{b}{c}\sin{\al c} & \frac{b}{c}(1-\cos{\al c})+\frac{a}{c}\sin{\al c} \\
        0 & \cos{\al c} & -\sin{\al c}  \\
        0 & \sin{\al c} & \cos{\al c}  \\
      \end{array}
    \right),
    \; c \neq 0
\right\},
\end{equation*}
\begin{equation*}\label{}
G_4=\left\{
e^A=\left(
      \begin{array}{ccc}
        1 & -\al b & \al a \\
        0 & 1 & 0  \\
        0 & 0 & 1  \\
      \end{array}
    \right),\;
    c = 0
    \right\}.
\end{equation*}
Both the cases for $\F_4$ can be generalized as it is shown in Table~1.

By similiar considerations we obtain the results in Table~1 for the classes $\F_9$ and $\F_{10}$.
%
%
%
%
%
%
%
%
%
%

\subsubsection{Lie groups as manifolds from the class $\F_8$}
According to Theorem A, we have
\begin{equation*}\label{}
[E_0,E_1]=\al E_2,\quad [E_0,E_2]=\al E_1,\quad [E_1,E_2]=0,
\end{equation*}
where $\al=-\lm$.

We get the matrix representation of $A$ for the considered Lie algebra $\g_8$, given in Table~1.

In this class $\tr{A}=0$ and we have three cases according to the sign of $\tr{A^2}=2\al^2\Dt$,
where $\Dt=2a^2-2b^2+c^2$. They can be generalized as it is shown in Table~1.

\section{Equipping of known Lie groups with almost contact B-metric structure}\label{sect-ex}

In this final section we equip some known Lie groups with almost contact B-metric structures of different types.
These considerations use results given in \cite{HM-Bia}.

\subsection{Examples 1, 2, 3}
In \cite{Kow}, there are considered the matrix Lie groups $G_\mathrm{I}$, $G_{\mathrm{II}}$ and $G_{\mathrm{III}}$ of the following form
\begin{equation*}\label{}
G_\mathrm{I}=
\left(
      \begin{array}{ccc}
        e^{-z}  & 0 & x \\
        0 & e^{z} & y \\
        0 & 0 & 1 \\
      \end{array}
    \right),\;\;
G_\mathrm{II}=
\left(
      \begin{array}{ccc}
        \cos{z}  & -\sin{z} & x \\
        \sin{z} & \cos{z} & y \\
        0 & 0 & 1 \\
      \end{array}
    \right),\;\;
G_\mathrm{III}=
\left(
      \begin{array}{ccc}
        1  & x & y \\
        0 & 1 & z \\
        0 & 0 & 1 \\
      \end{array}
    \right),
\end{equation*}
where $x,y,z\in\R$. These matrix groups and their automorphisms represent the only three types of generalized affine symmetric spaces denoted as $\mathrm{I}$, $\mathrm{II}$, $\mathrm{III}$.
The Lie group $G_\mathrm{I}$ is the group of hyperbolic motions of the plane $\R^2$,
$G_{\mathrm{II}}$ is the group of isometries of the plane $\R^2$ and $G_{\mathrm{III}}$ is the Heisenberg group.
The corresponding Lie algebras are determined by the commutators as follows:
\begin{equation*}\label{}
\begin{array}{llll}
\g_{\mathrm{I}}:\; &[X_1,X_3]=X_1,\quad &[X_2,X_3]=-X_2,\quad &[X_1,X_2]=0,\\
\g_{\mathrm{II}}:\; &[X_1,X_3]=-X_2,\quad &[X_2,X_3]=X_1,\quad &[X_1,X_2]=0,\\
\g_{\mathrm{III}}:\; &[X_1,X_3]=X_2,\quad &[X_2,X_3]=0,\quad &[X_1,X_2]=0.
\end{array}
\end{equation*}

The types of the Lie algebras $\g_{\mathrm{I}}$, $\g_{\mathrm{II}}$, $\g_{\mathrm{III}}$ according to the well-known
Bianchi classification in \cite{Bia1,Bia2} of three-dimensional real Lie algebras are
$\Bia{5}$, $\mathrm{Bia}(\mathrm{VII}_0)$, $\Bia{2}$, respectively.

Bearing in mind Table~1 for $\F_9$ and substituting $X_1=E_1$, $X_2=E_2$, $X_3=-E_0$, $\al=1$, we have that $G_\mathrm{I}$ can be considered as an almost contact B-metric manifold of the class $\F_9$ when the derived algebra is on $\Ker(\eta)$, otherwise the manifold belongs to $\F_1\oplus\F_{11}$. \cite{HM-Bia}

Considering Table~1 for $\F_4$ and substituting $X_1=E_1$, $X_2=E_2$, $X_3=E_0$, $\al=1$, we have that $G_\mathrm{II}$ can be considered as an almost contact B-metric manifold of the class $\F_4$ when the derived algebra is on $\Ker(\eta)$, otherwise the manifold belongs to $\F_4\oplus\F_{8}$ or $\F_4\oplus\F_{8}\oplus\F_{10}$. \cite{HM-Bia}

Depending on the way of equipping with an almost contact B-metric structure we can obtain $G_\mathrm{III}$ as a manifold in $\F_4\oplus\F_{10}$ or $\F_8\oplus\F_{10}$ when the derived algebra is on $\Ker(\eta)$ or $\Span(\xi)$, respectively. \cite{HM-Bia}

\subsection{Example 4}
Let us consider the well-known Lie group $SO(3)$ which is called the rotation group.
The matrix representation $A$ of the cor\-re\-spon\-ding Lie algebra has the form
\begin{equation*}\label{}
A=
\left(
      \begin{array}{ccc}
        0  & -\al z & \al y \\
        \al z & 0 & -\al x \\
        -\al y & \al x & 0 \\
      \end{array}
    \right),
\end{equation*}
where $\al$ is the angle of rotation and $x,y,z \in \R$.

In terms of the matrix exponential and according to the Rodrigues rotational formula,  the Lie group $SO(3)$ consists the matrices of the following form
\begin{equation*}\label{}
e^A=E+\sin{\al}\ A+(1-\cos{\al})\ A^2.
\end{equation*}

Then, the type of the considered Lie algebra is $\Bia{9}$, according to the Bianchi classification.
Hence, bearing in mind \cite{HM-Bia}, this Lie group can be considered as an almost contact B-metric manifold
of the class $\F_4\oplus\F_8\oplus\F_{10}$.



\begin{thebibliography}{33}

\bibitem{Bia1}
Bianchi, L., Sugli spazi a tre dimensioni che ammettono un gruppo
continuo di movimenti,
Memorie di Matematica e di Fisica della Societa
Italiana delle Scienze, Serie Terza, 11 (1898), 267-352.

\bibitem{Bia2}
Bianchi, L., On the three-dimensional spaces which admit a continuous group of motions,
Gen. Rel. Grav., 33 (2001), 2171-2253.

\bibitem{Blair}
Blair, D.E., Riemannian Geometry of Contact and Symplectic
Manifolds, Progress in Mathematics, Birkhauser, Boston, 2002.

\bibitem{GaBo}
Ganchev, G., Borisov, A., Note on the almost complex
manifolds with a Norden metric, C. R. Acad. Bulg. Sci., 39 (1986), 31-34.

\bibitem{GaMiGr}
Ganchev, G., Mihova, V. and Gribachev, K., Almost contact
manifolds with B-metric, Math. Balkanica (N.S.), 7 (1993) no. 3-4, 261-276.

\bibitem{GrMeDj}
Gribachev, K., Mekerov, D., Djelepov, G., Generalized
B-manifolds, C. R. Acad. Bulg. Sci., 38 (1985), 299-302.
(1985)



\bibitem{HM}
Manev, H., On the structure tensors of almost contact
B-metric manifolds, arXiv:\allowbreak{}1405.3088.

\bibitem{HM-Bia}
Manev, H., Almost contact B-metric structures and the Bianchi classification of the three-dimensional Lie algebras, arXiv:\allowbreak{}1412.8142.

\bibitem{IcoMek}
Manev, H. and Mekerov, D., Lie groups as 3-dimensional almost contact B-metric manifolds,
J. Geom., Online September 2014, DOI:10.1007/s00022-014-0244-0.

\bibitem{Gil}
Gilmore, R., Lie groups, Lie algebras and some of their applications, A Wiley-Interscience Publication, New York, 1974.


\bibitem{Kow}
Kowalski, O.,
Generalized Symmetric Spaces, Lecture Notes in Mathematics, 805 (Eds: A. Dold, B. Eckmann), Springer, Heidelberg, 1980.



\end{thebibliography}
\end{document}